\documentclass[11pt]{article}

\usepackage{lscape,tikz}

\definecolor{darkblue}{rgb}{0.2,0.2,0.71}
\usepackage{todonotes}
\usepackage{float}
\usepackage{colortbl}
\usepackage{multicol}
\usepackage{wrapfig}
\usepackage{color}

\usepackage{epsfig,cancel,ulem, amsthm,amssymb,amsmath}
\usepackage{color,soul}
\usepackage{graphicx}
\usepackage{yhmath}
\usepackage{mathdots}
\usepackage{MnSymbol}

\usetikzlibrary{decorations.markings,arrows}

\usepackage{wrapfig}
\usepackage{color}
\usepackage{graphicx,framed}
\usepackage[colorlinks=true, pdfstartview=FitV, linkcolor=darkblue, citecolor=darkblue, urlcolor=darkblue]{hyperref}
\definecolor{shadecolor}{rgb}{0.95, 0.95, 0.86}
\definecolor{darkgreen}{rgb}{0.2, 0.5,  0}

\def\&{\vspace{-5pt}&}

\def \C{\mathbb C}

\def \eqref#1{(\ref{#1})}
\def \& {&\hspace{-10pt}}

\newcommand{\bt}{\beta}

\renewcommand{\o}{\mathcal O}

\newcommand{\rf}{\rho_{ref} }  
\renewcommand{\O}{\Omega}

\newcommand{\br}{{\mathbb R}}

      \newcommand{\Lie}{\mathrm{Lie}}

\newtheorem{theorem}{Theorem}[section]
\newtheorem{example}[theorem]{Example}
\newtheorem{exercise}[theorem]{Exercise}

\newtheorem{lemma}[theorem]{Lemma}
\newtheorem{remark}[theorem]{Remark}

\newtheorem{proposition}[theorem]{Proposition}
\newtheorem{corollary}[theorem]{Corollary}
\newtheorem{definition}[theorem]{Definition}

\def\bt{\begin{theorem}}
\def\et{\end{theorem}}
\def\bc{\begin{corollary}}
\def\ec{\end{corollary}}
\def\bx{\begin{example}}
\def\ex{\end{example}}
\def\bxr{\begin{exercise}\small}
\def\exr{\end{exercise}}
\def\bl{\begin{lemma}}
\def\el{\end{lemma}}
\def\bd{\begin{definition}}
\def\ed{\end{definition}}
\def\bp{\begin{proposition}}
\def\ep{\end{proposition}}

\def\br{\begin{remark}}
\def\er{\end{remark}}

\def\be{\begin{equation}}
\def\ee{\end{equation}}

\def\&{\hspace{-15pt}&}
\def\bea{\begin{eqnarray}}
\def\eea{\end{eqnarray}}

\def\L{\mathrm {Lie}}

\def\1{{\bf 1}}

\newcommand{\F}{\mathbb F}

\makeatletter
\@addtoreset{equation}{section}
\makeatother

\begin{document}
\title{Frobenius manifolds on  orbits spaces}

\maketitle
\begin{center}
Zainab Al-Maamari
\footnote{Sultan Qaboos University, Muscat, Oman, s100108@student.squ.edu.om},
Yassir Dinar
\footnote{Sultan Qaboos University, Muscat, Oman, dinar@squ.edu.om}.
\end{center}

\maketitle

\begin{abstract}
The orbits space of an irreducible linear representation of a finite group is a variety whose coordinate ring is the ring of invariant polynomials. Boris Dubrovin proved that  the orbits space of the standard reflection representation of an irreducible finite  Coxeter group $\mathcal W$  acquires a natural   polynomial Frobenius manifold structure. We apply Dubrovin's method  on  various  orbits spaces of linear representations of  finite groups. We find some of them has non or several natural Frobenius manifold structures. On the other hand, these Frobenius manifold structures include  rational and trivial  structures  which are not  known to be related to the invariant theory of finite groups. \end{abstract}  

{\small \noindent{\bf Mathematics Subject Classification (2020) }  53D45}

{\small \noindent{\bf Keywords:} Invariant rings, Frobenius manifold, representations of finite groups, flat pencil of metrics, quotient singularities, orbifolds.}
\maketitle
\tableofcontents

\section{Introduction}

Frobenius manifold is a  geometric realization introduced by B. Dubrovin for a potential  satisfying a system of partial differential equations known  as  Witten-Dijkgraaf-Verlinde-Verlinde (WDVV) equations which describes the  module space of two dimensional topological field theory. Remarkably, Frobenius manifolds are also  recognized  in many other fields in mathematics like  invariant theory, quantum cohomology, integrable systems  and singularity theory \cite{DuRev}.   Briefly,  a Frobenius algebra is a commutative associative algebra with identity $e$ and a nondegenerate bilinear form $\Pi$ compatible with the product, i.e., $\Pi(a\circ b,c)=\Pi(a,b\circ c)$.
A Frobenius manifold is a manifold  with
a smooth structure of a Frobenius algebra on the tangent space  at any point  with certain compatibility conditions. Globally, we  require the metric $\Pi$ to be flat and the identity vector field $e$ is constant  with respect to its Levi-Civita connection. In this article, we  show that orbits spaces of some non-reflection representations of finite groups acquire Frobenius manifold structures.

We use the following notations and facts for a finite group $G$ and a  linear representation $\psi:G\to GL(V)$, where $V$ is a complex vector space. We denote by ${\mathbb C}[V]$ the ring of polynomial functions on $V$,   ${\mathbb C}[\psi]$  the subring of invariant polynomials in ${\mathbb C}[V]$, and ${\mathcal O} (\psi)$  the orbits space of the action of $G$ on $V$. Then  $ {\mathbb C}[\psi]$ is finitely generated by homogeneous polynomials and ${\mathcal O}(\psi)$ is a variety whose coordinate ring is  ${\mathbb C}[\psi]$ (\cite{Neusel}, \cite{dresken}). By Chevalley–Shephard–Todd theorem,  $ {\mathbb C}[\psi]$  is a polynomial ring if and only if $\psi$ is generated by pseudo-reflections. Let $(x^1,\ldots,x^n)$ be linear coordinates on $V$ and $f\in {\mathbb C}[\psi]$. Then the Hessian  $\mathrm H(f):={\partial^2 f\over \partial x^i \partial x^j}$ defines a bilinear from on the tangent space of ${\mathcal O} (\psi)$  and if $\det (\mathrm H(f))\neq 0$ then $f$ is a minimal degree invariant polynomial (\cite{Orlik}, page 6).   In this article, we will drop the word pseudo as all representations  will be   representations over complex vector spaces.  

Let ${\mathcal W}$ be a finite irreducible Coxeter group or  Shephard group of rank $r$ and ${\rho_{ref}}$ is the standard reflection representation of ${\mathcal W}$. Boris Dubrovin proved that the orbits space ${\mathcal O} ({\rho_{ref}})$ acquires  a polynomial Frobenius manifold  structure (\cite{DCG},\cite{DubAlmost}). This result led to the classification of irreducible semisimple polynomial Frobenius manifolds with positive degrees (see section \ref{first section} for more details). His method was used in \cite{polyZuo} when ${\mathcal W}$ is a Coxeter group of type $B_r$ or $D_r$ to construct $r$  Frobenius manifolds on ${\mathcal O} ({\rho_{ref}})$. In this article, we show that linear representations of finite groups are a valuable source to construct examples of Frobenius manifolds even if the representations are not  reflection representations. 


We mention that Dubrovin and his collaborators constructed Frobenius manifolds using invariant rings of infinite discrete groups being extensions of affine Weyl groups (\cite{Zang}, \cite{Stra}, \cite{Zuo2019}). However, we focus in this article  on linear representations of finite groups.

Let us fix a  finite group $G$ and a linear representation $\psi:G\to GL(V)$ of rank $r$. Then we summarize \textbf{Dubrovin's method} to construct Frobenius manifold structure on ${\mathcal O} (\psi)$ as follows:
\begin{enumerate}
    \item Fix homogeneous  invariants polynomial $f_1$ of  the minimal  degree.
    \item  Verify  that 
          the inverse of the Hessian  $H(f_1)$   defines a contravariant flat metric $\O_2$
    on some open subset $U$ of ${\mathcal O} (\psi)$. For example, this happens if  $\psi$ is a real representation (in this case degree $f_1$ equals 2) \cite{Fulton} or $\psi$ is the standard reflection representation of a  Shephard group \cite{Orlik}.    
    \item Construct another contravariant metric $\Omega_1$ which forms with $\O_2$ a regular quasihomogenius flat pencil of metrics (regular QFPM) on $U$ (see section \ref{preli2} for details).
    
    \item Then using a theorem   due to Dubrovin (see Theorem \ref{dub flat pencil} below), we get  a  Frobenius manifold structure on $U$ which depends on the representation $\psi$ of $G$ or ${\mathbb C}[\psi]$.
\end{enumerate}

\begin{definition} 
 By abuse of language, a Frobenius manifold structure obtained using Dubrovin's method will be called  a natural   Frobenius manifold structure on the orbits  space.
\end{definition}

Note that for a fixed metric $\O_2$, the problem of finding another metric $\O_1$ such that $(\O_2,\O_1)$ form a flat pencil of metric is not straightforward. For example, see the discussion on the classification of flat pencils of metrics related to the theory of Frobenius manifolds given in  \cite{Fordy}. We also observe that  an orbits space can have several natural Frobenius manifold structures. In this article, we will prove  the orbits spaces of the following representations posses natural Frobenius manifold structures: 

\begin{enumerate} 

\item The standard reflection representation  of a finite irreducible Coxeter group: We prove there is a natural  rational Frobenius manifold structure  different from the ones constructed in \cite{DCG} and \cite{polyZuo}. We give  details in section \ref{first section}.

\item  The non-standard irreducible representation of dimension  $r$ of  a Coxeter group of type $A_r$:  We show that it is a non-reflection representation and  we  construct certain $r$ algebraically independent invariant polynomials. Then,  we show that the orbits space   carries   natural rational  Frobenius manifold structures. We give the details in section  \ref{newrep}. 


    \item Irreducible representations of dihedral groups and dicyclic groups: These groups have only rank 1 and 2 irreducible representations. We will prove  that any rank 2  representation acquires  two natural Frobenius manifold structures.  See  section \ref{dicyclic and dihedral} for details.

    \item All  finite subgroups of the special linear group $SL_2({\mathbb C})$:  We get natural  polynomial and rational Frobenius manifold structures related to representations of  the dihedral groups. See section \ref{sl2}.  

    \item All finite subgroups of the special linear group $SL_3({\mathbb C})$ where the invariant rings are complete intersection: Dubrovin's method fails on some of them and we find natural trivial Frobenius manifold structures on others. We give details in section \ref{finite rank 3}. 

\end{enumerate}

As a consequence of this work, we noticed that  Frobenius manifold structures on orbits spaces of some non-reflection representations  appear in pairs.  Analyzing such pairs led us to the notion of the conjugate Frobenius manifold structures and we wrote the details on a separated article \cite{YZ1st}.  We  review this notion in section \ref{sec 2.3} and we show that the  conjugate of a natural Frobenius manifold structure is a natural Frobenius manifold structure.  

To make the article as self-contained as possible, we review in  section \ref{preli1} and \ref{preli2} the definition of Frobenius manifold and its relation with  flat pencils of metrics. 

\section{Flat pencil of metrics and Frobenius manifolds}
We review in this section the relation between Frobenius manifolds and flat pencil of metrics.
\subsection{Frobenius manifolds} \label{preli1}

Let $M$ be a Frobenius manifold with flat metric $\Pi$ and identity vector field $e$. In flat  coordinates $(t^1,...,t^r)$ for $\Pi$ where $e= \partial_{t^{r}}$ the compatibility conditions imply  that there exists a function $\mathbb{F}(t^1,...,t^r)$ which encodes  the Frobenius structure, i.e., the flat metric is given by
\begin{equation} \label{flat metric} \Pi_{ij}(t)=\Pi(\partial_{t^i},\partial_{t^j})=  \partial_{t^{r}}
\partial_{t^i}
\partial_{t^j} \mathbb{F}(t)\end{equation}
and, setting $\Omega_1(t)$ to be the inverse of the matrix $\Pi(t)$, the structure constants of the Frobenius algebra are given by 
\[ C_{ij}^k(t)=\Omega_1^{kp}(t)  \partial_{t^p}\partial_{t^i}\partial_{t^j} \mathbb{F}(t).\]
 Here, and in what follows, summation with respect to repeated
upper and lower indices is assumed. In this article, we assume the quasihomogeneity  condition for $\mathbb{F}(t)$ takes the form
\begin{equation}\label{quasihomog}
 E\mathbb F(t)=d_i t^i \partial_{t^i} \mathbb{F}(t) = \left(3-d \right) \mathbb{F}(t);~~d_{r}=1.
\end{equation}
 The vector field $E= d_i t^i \partial_{t_i}$ is known as  Euler vector field and it  defines  the degrees $d_i$ and  the charge $d$ of  $M$. 
The  associativity  of the Frobenius
algebra implies that the potential $\mathbb{F}(t)$ satisfies WDVV  equations, i.e.,
\begin{equation} \label{frob}
 \partial_{t^i}
\partial_{t^j}
\partial_{t^k} \mathbb{F}(t)~ \Omega_1^{kp} ~\partial_{t^p}
\partial_{t^q}
\partial_{t^n} \mathbb{F}(t) = \partial_{t^n}
\partial_{t^j}
\partial_{t^k} \mathbb{F}(t) ~\Omega_1^{kp}~\partial_{t^p}
\partial_{t^q}
\partial_{t^i} \mathbb{F}(t),~~ \forall i,j,q,n. 
  \end{equation}
We say $M$ is a polynomial (resp. rational) if $\mathbb{F}(t)$ is a polynomial (resp. rational) function.
\begin{definition}
Let $M$ and $\widetilde M$ be two Frobenius manifolds  with flat metrics $\Pi$ and $\widetilde \Pi$. Let  $\mathbb{F}$ and $\widetilde{\mathbb{F}}$ be the corresponding potentials, respectively. We say $M$ and $\widetilde{M}$ are (locally) equivalent if there are open sets $U\subseteq M$ and $\widetilde U \subseteq \widetilde M $ with a local diffeomorphism $\phi:U\to \widetilde U$ such that \begin{equation}\label{equivalent} \phi^* \widetilde\Pi=c^2 \Pi,\end{equation} for some nonzero constant $c$, and $\phi_*:T_tU\to T_{\phi(t)} \widetilde U$ is an isomorphism of Frobenius algebras. 
\end{definition} 
Note that, if $M$ and $\widetilde M$, are equivalent Frobenius structures then it is not necessary that $\phi^* \widetilde{\mathbb F}=\mathbb F$ \cite{DCG}. 

\subsection{Flat pencil of metrics}\label{preli2}

We  review the relation between flat pencils of metrics and Frobenius manifolds outlined  in \cite{Du98}.

Let $M$ be a smooth manifold of dimension $r$ and  fix  local coordinates $(u^1, . . . , u^r)$ on $M$. 

\bd \label{contra metric} A symmetric bilinear form $(. ,. )$ on $T^*M$ is called a contravariant
metric if it is invertible on an open dense subset $M_0 \subseteq M$. We define the contravariant Christoffel symbols $\Gamma^{ij}_k$  for a contravariant
metric $(. ,. )$ by
\[
\Gamma^{ij}_k:=-\Omega^{im} \Gamma_{mk}^j
\]
where $\Gamma_{mk}^j$ are the  Christoffel symbols of the metric $<. ,. >$ defined on $TM_0$ by the inverse of the matrix $\Omega^{ij}(u)=(du^i, du^j)$.
We say the metric $(.,.)$ is flat if  $<. ,. >$ is flat.
\ed

Let $(. ,. )$ be a contraviariant metric on $M$ and set $\O^{ij}(u)=(du^i, du^j)$. Then we will use $\Omega$ to refer to the metric and $\Omega(u)$ to refer to its  matrix  in the coordinates. In particular, the Lie derivative of $(. ,. )$ along a vector field $X$  will be written $\Lie_X \Omega$ while $X\Omega^{ij}$ means the vector field $X$ acting on the entry $\Omega^{ij}$. The Christoffel symbols given in Definition \ref{contra metric} determine for $\O$ the contravariant (resp. covariant) derivative  $\nabla^{i}$ (resp.  $\nabla_{i}$) along the covector $du^i$ (resp. the vector field $\partial_{u^i}$). They are related by the identity  $\nabla^{i}=\O^{ij}(u) \nabla_{j}$.

\bd 
A flat pencil of metrics (FPM)  on $M$ is a pair $(\Omega_2,\Omega_1)$ of 
 two flat contravariant metrics $\O_2$ and $\O_1$ on $M$ satisfying 
 \begin{enumerate}
     \item $\O_2+\lambda \O_1$ defines a flat metric on $T^*M$ for  a generic constant $\lambda$,
     \item the Christoffel symbols of $\O_2+\lambda \O_1$ are   $\Gamma_{2k}^{ij}+\lambda \Gamma_{1k}^{ij}$, where   $\Gamma_{2k}^{ij}$ and $ \Gamma_{1k}^{ij}$ are the Christoffel symbols of $\O_2$ and $\O_1$,  respectively. 
 \end{enumerate}  
\ed  
\bd \label{FPM}  A  flat pencil of metrics $(\O_2,\O_1)$ on  $M$ is called quasihomogeneous flat pencil of metrics (QFPM) of  degree $d$ if there exists a function $\tau$ on $M$ such that the 
 vector fields $E$ and $e$ defined by 
\begin{eqnarray} \label{tau flat pencil} E&=& \nabla_2 \tau, ~~E^i
=\O_2^{ij}(u)\partial_{u^j}\tau
\\\nonumber  e&=&\nabla_1 \tau, ~~e^i
= \O_1^{ij}(u)\partial_{u^j}\tau  \end{eqnarray} satisfy 
\be \label{vector fields} [e,E]=e,~~ \Lie_E \O_2 =(d-1) \O_2,~~ \Lie_e \O_2 =
\O_1~~\mathrm{and}~~ \Lie_e\O_1
=0.
\ee 
Such a QFPM is  \textbf{regular} if  the
(1,1)-tensor
\begin{equation}\label{regcond}
  R_i^j = \frac{d-1}{2}\delta_i^j + {\nabla_1}_i
E^j
\end{equation}
is  nondegenerate on $M$.
\ed 

We will use the following source for FPM.

\begin{lemma}\label{almost linear} \cite{DCG}  Let $\Omega_2$ be a contravariant flat metric on $M$. Assume that in the coordinates $(u^1,...,u^r)$,  $\Omega^{ij}_2(u)$ and   $\Gamma^{ij}_{2k}(u)$ depend almost linearly on $u^{r}$. Suppose that $\Omega_1:=\Lie_{\partial_{u^r}} \Omega_2=\partial_{u^r} \Omega_2(u)$ is nondegenerate. Then 
 $(\Omega_2,\Omega_1)$ form a FPM. The Christofell symbols of $\Omega_1$ has the form $
 \Gamma^{ij}_{1k}(u)=\partial_{u^r} \Gamma^{ij}_{2k}(u).
$
\end{lemma} 

If $M$ is a Frobenius manifold then $M$ has a QFPM of degree $d$ but it does not necessarily satisfy the regularity condition \eqref{regcond} \cite{Du98}. In the notations of section \ref{preli1}, the QFPM consists of the intersection form  ${\Omega}_2(t)$ and the flat metric ${\Omega}_1(t)$ where  \begin{equation} \label{intesection form}
{\Omega}_2^{ij}(t):=(d-1+d_i+d_j)\Omega^{i\alpha}_1\Omega^{j\beta}_1
\partial_{t^\alpha}
\partial_{t^\beta} \mathbb{F}.
\end{equation}
 Furthermore,   $\tau=\Pi_{i1}t^i$ and $E$ with $e$ are defined by \eqref{tau flat pencil} and satisfy equations \eqref{vector fields} . The converse is given by the following theorem

\bt\cite{Du98}\label{dub flat pencil} Let $M$ be a manifold carrying a regular QFPM  $(\Omega_2,\Omega_1)$ of degree $d$. Then there exists a unique Frobenius manifold structure on $M$ of charge $d$ where $(\Omega_2,\Omega_1)$ is the associated QFPM. 
\et

\section{Conjugate Frobenius Manifold and Dubrovin's method}\label{sec 2.3}

We begin this section with a theorem proved in \cite{YZ1st} which leads to the notion of conjugate Frobenius manifold structure. Then we will prove that the conjugate natural Frobenius manifold structure constructed on an orbits spaces is also natural.

\bt \cite{YZ1st} \label{dual Frob manif}
Let $M$ be a Frobenius manifold with the Euler vector field $E$ and the identity vector field $e$. Suppose the associated QFPM  is regular of degree $d$ with  a function $\tau$. Assume that
\be \label{new cond}e(\tau)=0 \ \ \text{and} \ \ E(\tau)=(1-d)\tau. 
\ee Then  we can construct another  Frobenius manifold structure on $M\backslash \{\tau=0\}$ of degree $2-d$. Moreover, we can apply the same  method to the new Frobenius manifold structure  and it leads to the original Frobenius manifold structure. 
\et  

For a fixed Frobenius manifold  the new structure that can be obtained using Theorem \ref{dual Frob manif} will be called the conjugate Frobenius manifold structure.

Let $M$ be a Frobenius manifold of degree $d$.  Let $T=(\O_2,\O_1)$ be the associated QFPM  with a function $\tau$, the Euler vector field $E$ and the identity vector field $e$. Suppose it satisfies the hypothesis of Theorem \ref{dual Frob manif}. Then the QFPM associated to the conjugate Frobenius manifold structure has the form $\widetilde T:=(\O_2,\widetilde \O_1)$ where $\widetilde \O_1:=\Lie_{\widetilde e}\O_2 $ and the vector field $\widetilde e:= \tau^\frac{2}{1-d} e$ \cite{YZ1st}.

Let us adapt the notations of section \ref{preli1} and  assume  $\Pi_{ij}=\delta_{i+j}^{r+1}$, i.e., the potential  $\F$ has the standard form 
\be \label{norm potential}
 \F(t) = \frac{1}{2} (t^r)^2 t^1 + \frac{1}{2} t^r \sum_{i=2}^{r-1} t^i t^{r-i+1} + G(t^1,...,t^{r-1}).
\ee 
Then we get the following consequence of Theorem \ref{dual Frob manif}.

\bt \cite{YZ1st} \label{main thm}
Let $M$ be a Frobenius manifold with charge $d\neq 1$. Suppose in the flat coordinates $(t^1,\ldots,t^r)$, the potential $\mathbb F(t)$ has the standard form \eqref{norm potential} and the quasihomogeneity condition takes the form  \eqref{quasihomog} with $d_i\neq \dfrac{d_1}{2}$ for every $i$.  Then we can construct the  conjugate   Frobenius manifold structure on $M\backslash\{t^1=0\}$.    Moreover,  flat coordinates  for the conjugate Frobenius manifold  are 
\be \label{change coord}
s^1= -t^1 , \ \ s^i= t^i (t^1)^{\frac{d_1-2d_i}{d_1}}\ \ for \ \ 1<i< r, \ \ s^r= \frac{1}{2} \sum_{i=1}^{r}  t^i t^{r-i+1} (t^1)^{\frac{-2}{d_1}-1}.
\ee 
In addition, the corresponding potential   equals the potential  obtained by applying the inversion symmetry to $\mathbb F(t)$  and it is given by   \be\label{F in 3 coord}
\widetilde{\F}(s) = (t^1)^{ \frac{-4}{d_1}}  \left(\F(t^1,\ldots,t^r) -\frac{1}{2} t^r \sum_{1}^{r} t^i t^{r-i+1}\right). 
\ee 
The degrees $\widetilde d_i$ and the charge $\widetilde d$  of the conjugate Frobenius manifold structure are given by 
\be \label{degrees inv} \widetilde{d}_1=-d_1,\ \ \widetilde{d}_r=1, \ \ \widetilde{d}_i = d_i-d_1 \ \ for \ \ 1< i < r,~\widetilde d=2-d.\ee
\et 

See (\cite{DuRev}, Appendix B) for details about inversion symmetry of solutions to WDVV equations.  Form the point of view of this  article, Theorem \ref{main thm} explains the appearance  of pairs of natural Frobenius manifold structures on orbits space of some linear representations of finite groups.

\bt \label{conj on orbits}
Let $M$ be the orbits space of a linear representation  of a finite group. Assume $M$ inherits a natural  Frobenius manifold structure which  has a conjugate Frobenius manifold structure. Then the conjugate Frobenius manifold structure on $M$ is also natural.  
\et 

\begin{proof}
Let $T=(\O_2,\O_1)$ be the associated QFPM of the Frobenius manifold structure on $M$ which is obtained using Dubrovin's method. Then  $\O_2$ is defined using  the Hessian of a minimal invariant polynomial $f_1$. The QFPM  associated to the conjugate Frobenius manifold  has the same intersection  form $\O_2$ and hence it constructed by Dubrovin's method. 
\end{proof}
 

 For convenience, we write in  examples, indices of coordinates using subscripts instead of superscripts.

\begin{example}\label{Triv1}
The potential \begin{equation}\label{triv1}
\F=\frac{t_1^3}{6}-\frac{1}{2} t_2^2 t_1+\frac{1}{2} t_2^2 t_3+\frac{1}{2} t_1 t_3^2.
\end{equation}
defines two inequivalent trivial Frobenius manifold structures, i.e., both have   charge $d=0$ and   Euler vector field $E=\sum t_i \partial_{t_i}$. Setting the identity vector field to be  $\widehat e= \partial_{t_1}$, $\mathbb F$ defines a Frobenius manifold structure $\widehat T_3$ whose associated regular QFPM does not satisfy condition \eqref{new cond}, i.e., it does not have a conjugate structure.  While fixing the identity vector field $e=\partial_{t_3}$, we get a  Frobenius manifold structure $T_3$ which has conjugate. The associated regular QFPM $(\O_2,\O_1)$ has $\O_1^{ij}(t)=\delta^{i+j}_3$ while 
\[
\Omega_2(t) =\left(
\begin{array}{ccc}
 t_1 & t_2 & t_3 \\
 t_2 & t_3-t_1 & -t_2 \\
 t_3 & -t_2 & t_1 \\
\end{array}
\right)
\]
Setting 
\[
s_1=-t_1,\ \ s_2=\frac{t_2}{t_1},\ \  s_3=\frac{t_2^2}{2 t_1^3}+\frac{t_3}{t_1^2}
\]
the conjugate QFPM has $\widetilde \O_1^{ij}(s)=\delta^{i+j}_3$ and 
\[\O_2(s)=\left(
\begin{array}{ccc}
 -s_1 & 0 & s_3 \\
 0 & s_3+\frac{3 s_2^2}{2 s_1}+\frac{1}{s_1} & -\frac{s_2^3}{s_1^2}-\frac{2 s_2}{s_1^2} \\
 s_3 & -\frac{s_2^3}{s_1^2}-\frac{2 s_2}{s_1^2} & \frac{3 s_2^4}{4 s_1^3}+\frac{3 s_2^2}{s_1^3}-\frac{1}{s_1^3} \\
\end{array}
\right)\]
The potential of the  conjugate 
Frobenius manifold structure reads
\begin{equation}\label{triv1conj}
\widetilde{\F}(s)=\frac{-1}{6 s_1} +\frac{s_2^2}{2 s_1}  +\frac{s_2^4}{8 s_1}+\frac{1}{2} s_2^2 s_3+\frac{1}{2} s_1 s_3^2.\end{equation}
Here  ${\widetilde{E}}=-s_1 \partial_{s_1}+s_3 \partial_{s_3}$ and  $\widetilde{E} \widetilde{\F}= \widetilde{\F}$. 
\end{example}

\section{Coxeter groups} 
\subsection{The standard   reflection representation} \label{first section}

In  this section, we  recall  the standard reflection representations of irreducible finite Coxeter groups and review the construction of natural  Frobenius manifolds on their orbits space. Then we classify  those   having conjugate Frobenius manifold structures. Note that the conjugate Frobenius manifold structures will be rational and they are not known to be related to invariant theory of finite groups. 

We fix  an irreducible finite Coxeter system $(\mathcal W,S)$ of rank $r$, i.e.,  
\begin{equation}  \mathcal W=<S|\, (s s')^{m(s,s')} =1;\, \forall s,s'\in S>,~~ r=|S|. \end{equation}

Let $V$ be the formal vector space over ${\mathbb C}$ with  basis $\{ \alpha_s \ | \ s\in S \}$. Then the standard reflection representation of $\mathcal W$ is defined by
\begin{align*}
{\rho_{ref}} &: \mathcal W \rightarrow GL(V),~~ s\mapsto R_s,~ s\in S.\\ R_s(v)&:= v -2B(\alpha_s, v) \alpha_s,~~v \in V,~~B(\alpha_s, \alpha_{s'} ) := -\cos \frac{\pi}{m(s,s')}.
\end{align*} Here $B$ is the standard positive-definite Hermitian  form  on $V$ which is  invariant under $\rho_{ref}$.  By Chevalley–Shephard–Todd theorem,  the invariant ring ${\mathbb C}[\rho_{ref}]$ is  a polynomial ring generated by $r$   homogeneous polynomial.  We fix  generators  $u^1,...,u^r$ for ${\mathbb C}[\rho_{ref}] $. We assume $\deg u^i =\eta_i$ and  
\begin{equation} 2=\eta_1 < \eta_2 \leq \eta_3\leq \ldots \leq\eta_{r+1}< \eta_r.\end{equation} These degrees are  uniquely determined by the group $\mathcal W$ \cite{Hum}.

We  assume   $u^1$ equals the quadratic from of $B$. Hence, the inverse of  the Hessian of $u^1$  defines a flat contravariant  metric $\Omega_2$ on ${\mathcal O}(\rho_{ref})$. It is easy to prove that $\Omega_2(u)$ is almost linear in $u^r$ by analysing the degrees of $\Omega^{ij}_2(u)$.  We fix the vector field $e=\partial_{u^r}$. Note that changing the generators of ${\mathbb C}[\rho_{ref}]$, $e$ is uniquely defined up to a constant factor. Setting $\Omega_1 :=\Lie_{e}\Omega_2$ Dubrovin  proved  that $T:=(\Omega_2,\O_1)$ is a regular QFPM of charge $\eta_r-2\over \eta_r$ \cite{Du98}. In this case,  $\tau={1\over \eta_r} u_1$ and  the vector field $E$  is given by $E={1\over \eta_r}\sum_i\eta_i u^i \partial_{u^i}$. This result initiated what we call Dubrovin's method. We observe that $E$ is uniquely defined and does not depend on the choice of invariants $u^i$.  Also, we   mention that the flat metric $\Omega_1$ was  studied by  K. Saito \cite{Saito}, \cite{Saito1} and his results was very important to the work \cite{DCG}.  We restate Dubrovin's theorem.

\begin{theorem}  (\cite{DCG},  \cite{Du98})\label{polyFrob}
 The FPM $(\Omega_2,\Omega_1)$ defines a unique (up to equivalence) natural polynomial Frobenius manifold on ${\mathcal O}(\rho_{ref})$ with  degrees $\eta_i\over \eta_r$ and charge $\eta_r-2\over \eta_r$.
\end{theorem}

The following theorem was conjectured by Dubrovin  and proved by C. Hertling. 

\begin{theorem}  \cite{Hert} \label{polyFrob2}
Any irreducible massive polynomial Frobenius manifold with  positive degrees is isomorphic to a polynomial Frobenius manifold constructed by Theorem \ref{polyFrob} on the orbit space of the standard reflection representation of an  irreducible finite Coxeter group.
\end{theorem}

The following theorem grantees the existence of  another natural  Frobenius manifold structure  on ${\mathcal O}(\rho_{ref})$.  

\begin{theorem} \label{dual  poly}
The polynomial Frobenius manifold constructed by Theorem \ref{polyFrob} on the orbits space  $\o(\rho_{ref})$ has a conjugate Frobenius manifold structure. Thus, we get a rational natural Frobenius manifold structure on ${\mathcal O}(\rho_{ref})$. \end{theorem}
\begin{proof}
 There exist invariant polynomials $t^1,\ldots, t^r$ which form flat coordinates and the potential has the form \eqref{norm potential} \cite{DCG}.  From the structure of the degrees, we can and we will apply Theorem \ref{main thm} to get a rational  conjugate  Frobenius manifold. The last statement is a consequence of Theorem \ref{conj on orbits}.
\end{proof}

Let us assume $\mathcal W$ is of type $B_r$. Then  Dafeng Zuo  obtained $r$ Frobenius manifold  structures on ${\mathcal O}(\rho_{ref})$ by fixing  certain generators $z^1,\ldots, z^r$ for ${\mathbb C}[\rho_{ref}]$ \cite{polyZuo}. Under these generators,  $\Omega_2(z)$ and its Christoffel symbols  $\Gamma^{ij}_{2k}(z)$ are almost linear in each $z^k$, $k=1,2,\ldots,r$. Then he proved  that Lemma \ref{almost linear} can be applied and  he constructed $r$ rational Frobenius manifold  structures using the flat pencils of metrics  $\widehat T_k:=(\Omega_2,\Lie_{\partial_{z^k}}\Omega_2)$. He also proved that the same Frobenius manifold structures can be constructed when $\mathcal W$  is of type  $D_r$. Even it is not written explicitly in \cite{polyZuo}, We confirm that they are natural Frobenius manifold structures as  each $\widehat T_k$ is   regular QFPM of degree $1-{1\over k}$ with $\tau={1\over 4k} z^1$. Here $e=\partial_{z^k}$. Thus, we can obtain these  Frobenius manifolds directly using Theorem \ref{dub flat pencil}.  Here the structure of  Zuo's theorem
\begin{theorem}\cite{polyZuo}\label{zuo}
There exists a unique natural  Frobenius structure for each $1 \leq k \leq r$ of charge $d=1-\frac{1}{k}$ on the orbit space ${\mathcal O}(\rho_{ref})$ when $\mathcal W$ is of type $B_r$ and $D_r$ polynomial in $t^1,t^2,\ldots,t^r,\frac{1}{t^r}$ such that:
\begin{enumerate}
    \item The identity vector field is  $e=\frac{\partial}{\partial z^k}=\frac{\partial}{\partial t^k}$.
     \item  The Euler vector field is $E=  \sum_{i=1}^{r} d_i  t^i \partial_{t^i} $, where
     \[ d_1=\frac{1}{k}, \ \ d_i=\frac{i}{k}\ \ for \  \  2\leq i\leq k, \ \ d_i=\frac{2k(r-i)+r}{2k(r-k)} \ \ for \  \ k+1\leq i\leq r. \]
  \item The assciated QFPM is $\widehat T_k$.
\end{enumerate}

\end{theorem}

Note that when $k=1$, $\widehat{T}_1$ does not satisfy condition \eqref{new cond}. Thus the corresponding Frobenius manifold structure has no conjugate. For $k>1$, we get the following theorem.

\begin{theorem}\label{dual zuo}
For $k>1$, the natural Frobenius manifold structure corresponding to $\widehat T_k$ constructed by Theorem \ref{zuo} has a conjugate Frobenius manifold structure which is also natural. 
\end{theorem} 
\begin{proof}
 Similar to the proof of Theorem \ref{dual  poly}, we apply Theorem \ref{main thm} and Theorem \ref{conj on orbits}.
 
\end{proof}

Considering Theorem \ref{polyFrob2}, let  $K$ be the type of $\mathcal W$,  then we say a  Frobenius manifold  is of type $K$ (rep. of type $\widetilde K$) if it isomorphic to a natural polynomial   Frobenius manifold (resp.  a natural conjugate Frobenius manifold) constructed on ${\mathcal O}(\rho_{ref})$  by Theorem \ref{polyFrob} (resp. Theorem \ref{dual  poly}).

\begin{example} \label{Dualtiy poly}
 We list in Table \ref{table2} all Frobenius structures constructed on ${\mathcal O}(\rho_{ref})$ when $\mathcal W$ is of rank 3 using the above theorems.  We borrow  the potentials of Frobenius structures of type $A_3$, $B_3$ and $H_3$ from \cite{Du98}. From these potentials, we find Frobenius manifold structures of type $\widetilde A_3$, $\widetilde B_3$ and $\widetilde H_3$ using the formula \eqref{F in 3 coord}. Then applying Theorem \ref{zuo} to a  Coxeter group of type $B_3$, we get a Frobenius manifold  of type $B_3$ (resp. $A_3$) when  $k=3$ (resp. $k=2$). For  $k=1$, we get a rational Frobenius manifold $B_3^1$  which has no conjugate.

  \setlength{\arrayrulewidth}{0,2mm}
\setlength{\tabcolsep}{2pt}
\renewcommand{\arraystretch}{0.5}
  \begin{table}[H]
\centering

\begin{tabular}{| c | c | c | c |}
\hline
Notations & ${\mathbb{F}}(t_1,t_2,t_3)$& $d_1,d_2,d_3$ & d  \\ [0.7ex]
\hline\hline
 $A_{3}$ & $\frac{1}{2} t_3^2 t_1  + \frac{1}{2} t_2^2 t_3 + \frac{1}{4} t_1^2 t_2^2+ \frac{1}{60} t_1^5$ &  $\frac{1}{2},\frac{3}{4},1$ & $\frac{1}{2}$\\[2ex] 

$\widetilde{A}_{3}$ & $ \frac{1}{2}t_3^2 t_1 +\frac{1}{2}  t_2^2 t_3+\frac{t_2^4}{8 t_1}+\frac{t_2^2}{4 t_1^2}-\frac{1}{60 t_1^3}$ &  $\frac{-1}{2},\frac{1}{4},1$ & $\frac{3}{2}$ \\ [1ex] 
\hline
 $B_{3}$ & $\frac{1}{2} t_3^2 t_1+\frac{1}{2} t_2^2
   t_3+\frac{1}{6} t_2^2 t_1^3+\frac{1}{6} t_2^3 t_1+\frac{1}{210} t_1^7$ & $\frac{1}{3},\frac{2}{3},1$ & $\frac{2}{3}$\\
[2ex] 
$\widetilde{B}_{3}$ & $\frac{1}{2}t_3^2 t_1+\frac{1}{2} t_2^2 t_3 
   + \frac{t_2^4}{8 t_1}+\frac{t_2^3}{6 t_1^2}-\frac{t_2^2}{6 t_1^3}-\frac{1}{210 t_1^5}$ & $\frac{-1}{3},\frac{1}{3},1$ & $\frac{4}{3}$\\[1ex] 
\hline
$H_{3}$ &$\frac{1}{2} t_3^2 t_1+\frac{1}{2}
   t_2^2 t_3+\frac{1}{20} t_2^2 t_1^5+\frac{1}{6} t_2^3 t_1^2+\frac{1}{3960} t_1^{11}$ &  $\frac{1}{5},\frac{3}{5},1$ & $\frac{4}{5}$ \\
[2ex] 
$\widetilde{H}_{3}$ & $\frac{1}{2} 
   t_3^2 t_1 +\frac{1}{2}t_2^2 t_3 +\frac{t_2^4}{8 t_1}-\frac{t_2^3}{6 t_1^3}-\frac{t_2^2}{20 t_1^5}-\frac{1}{3960 t_1^9}$ &  $\frac{-1}{5},\frac{2}{5},1$ & $\frac{6}{5}$ \\[1ex] 
\hline
$B_3^1$& $ \frac{1}{2}t_3^4+\frac{3}{2}
   t_1 t_2
   t_3+\frac{1}{8}t_1^3+{1\over 16}\frac{t
   _2^3}{ t_3}$ &  $1,\frac{3}{4},\frac{5}{4}$ & $0$ \\
\hline
\end{tabular} 

\caption{Frobenius manifolds on orbits spaces of reflection groups of rank 3} 
\label{table2}
\end{table}

\end{example}

\subsection{Sign times reflection representation}\label{newrep}

We keep the notations of the last section and we assume $\mathcal W$ is of type $A_r$.  We study an irreducible    representation $\rho_{new}$ of $\mathcal W$ which can be defined  using the sign representation and the representation $\rho_{ref}$. The definition  will enable us to construct $r$ invariant polynomials of   $\rho_{new}$. We will prove the invariant ring  ${\mathbb C}[\rho_{new}]$ is not  a polynomial ring when $r>2$. We recall that the degrees of a complete set of generators of $\C[\rho_{ref}]$ are $2,3,\ldots,r+1$.

 We consider the sign representation of $\mathcal W$,  $\rho_{sign} : \mathcal W \rightarrow {\mathbb C}^*$ defined by  sending each element $s\in S $ to $-1$.  Then we define the representation $\rho_{new}$ of $\mathcal W$ by 
 \begin{equation} \rho_{new}: \mathcal W \rightarrow GL({\mathbb C} \otimes V),~~ \rho_{new} (w)= \rho_{sign}(w) \otimes \rho_{ref}(w),~~\forall w\in \mathcal W .\end{equation}
  Note that $\rho_{new}$ is a real representation  of rank $r$. The following proposition  proves that $\rho_{new}$ is an irreducible representation. 

\begin{proposition}
The new representation $\rho_{new}$ is an irreducible representation of $\mathcal W$. Moreover,   $\rho_{new}$ and $\rho_{ref}$ are  isomorphic when  $r=2$ and different otherwise.   
\end{proposition}
\begin{proof}
Recall that if $\chi_\psi$ denotes the character of a representation $\psi$ of a finite group $G$, then $\psi$ is irreducible if and only if
\cite{BST}
\begin{equation} 
\frac{1}{|G|} \sum_{g \in G} \chi_{\psi}(g) \overline{\chi_{\psi}(g)}=1.\end{equation}
Note that $\rho_{ref}$ and $\rho_{sign}$ are irreducible representations and $$\chi_{\rho_{new}} (w) =\chi_{\rho_{sign}} (w) \chi_{\rho_{ref}} (w).$$ Then
\begin{align}
\frac{1}{|\mathcal W|} \sum_{w \in \mathcal W} \chi_{{\rho_{new}}}(w) \overline{\chi_{{\rho_{new}}}(w) } &=  \frac{1}{|\mathcal W|} \sum_{w \in \mathcal W} (\chi_{\rho_{sign}}(w)  \chi_{\rho_{ref}}(w)) \overline{(\chi_{\rho_{sign}}(w)  \chi_{\rho_{ref}}(w))}\\ \nonumber 
&= \frac{1}{|\mathcal W|} \sum_{w \in \mathcal W} (\chi_{\rho_{ref}}(w) \overline{\chi_{\rho_{ref}}(w)})=1. \nonumber 
\end{align}
 For the second part,  note that for any generator  $s\in S$, $\chi_{{\rho_{new}}}(s)=-\chi_{\rho_{ref}}(s)=-(r-2)$. Hence, the two representations are different when $r\neq 2$. For $r=2$, we can check that $\rho_{new}$ is equivalent to       $\rho_{ref}$ by direct computations. 

\end{proof}

For the remainder of this section we assume the rank $r>2$. Recall that the Coxeter group of type  $A_r$ is isomorphic to the symmetric group $S_{r+1}$. Thus, irreducible representations of $A_r$ are in one to one correspondence with the partition of $r+1$. For a given partition $\lambda$ of $r+1$, the corresponding  irreducible representation can be constructed using Young tableaux associated to $\lambda$ \cite{Fulton}. Under this construction, the reflection representation $\rho_{ref}$ is associated with the partition $[r,1]$, $\rho_{sign}$ is associated with the partition $[r+1]$ while $\rho_{new}$ is associated with $[2,1,1,\ldots,1]$. The character of each representation is given by Frobenius formula \cite{Fulton}. We use this formula to prove the following proposition.

\begin{proposition} \label{new not ref}
The irreducible representation $\rho_{new}$ is not a reflection representation. In particular the ring ${\mathbb C}[\rho_{new}]$ is not a polynomial ring. 
\end{proposition} 

\begin{proof}
Assume that $\rho_{new}$ is a reflection representation. Then, it is generated by a set of involutions $w_1,\ldots,w_r$. Since $\rho_{new}$ is a real representation, we must have  $\chi_{\rho_{new}}(w_i)=r-2$. From  $\rho_{new}(w_i)=\rho_{sign}(w_i) \rho_{ref} (w_i)$, we have  $\rho_{sign} (w_i)=-1$, since if  $\rho_{sign} (w_i)=1$,  then $\rho_{ref} (w_i)$ is a reflection and we get a contradiction. Thus,  $\chi_{\rho_{ref}}(w_i)=2-r$. In the one-to-one correspondence between conjugacy classes of $S_{r+1}$ and partitions of $r+1$, $w_i$  corresponds to a partition of the from $[2,2,\ldots,2,1,1\ldots,1]=[2^p,1^q]$ with $2p+q=r+1$, $p>0$.  Using Frobenius formula, $\chi_{\rho_{ref}}(w_i)$ equals  the coefficient of  $x^{r+1} y$ in the expansion $(x-y)(x^2+y^2)^p(x+y)^{r+1-2p}$.  Hence, $\chi_{\rho_{ref}}(w_i)=r-2p$. Using the fact that $2p\leq r+1$ and $\chi_{\rho_{ref}}(w_i)=2-r$ we get $r\leq 3$. However, the case $r=3$ is excluded by direct computations.  
\end{proof}

 We study the ring  ${\mathbb C}[\rho_{new}]$ in order to use Dubrovin's method. 
We fix a basis  $e_1,e_2, \ldots, e_r$ for $V$ and let $x^1,\ldots,x^r$  be the dual basis satisfying $x^i(e_j)=\delta^i_j$. Then $ \tilde{e_i}:=\textbf{1} \otimes e_i$, $i=1,...,r$ form a basis of ${\mathbb C} \otimes V$ and we get a natural isomorphism  
\begin{equation}\label{iso} \theta:{\mathbb C} \otimes V \rightarrow V, \, \tilde{e_i} \mapsto e_i. \end{equation} 
Then the pullback  $\tilde{x}^i=\theta^*(x^i)$ defines the dual basis of $\tilde{e}_i$. Let $w\in \mathcal W$ and $a_i^j$ be the matrix of $\rho_{ref}(w)$ under the basis $e_i$. Then  
 $ \rho_{new} (w)( \tilde{e_i})= \rho_{sign} (w) \textbf{1}\otimes \rho_{ref}(w)e_i  =
\rho_{sign} (w) a_i^j  \ \tilde{e_j}$. Therefore,  $\rho_{new}(w)=\rho_{sign}(w) \rho_{ref}(w)$. 

 \begin{lemma}\label{lem1}
Let $w\in \mathcal W$ with $
\rho_{sign}(w) \rho_{new}(w)\notin \rho_{ref}(\mathcal W)$ and $ f\in  {\mathbb C}[\rho_{ref}]$ be homogeneous polynomial.   Then  
\begin{equation} w\cdot \theta^*(f)= (\rho_{sign}(w))^{deg(f)} \theta^*(f).\end{equation}
In particular,  if degree $f$ is even then $\theta^*(f)\in {\mathbb C}[\rho_{new}]$.
 \end{lemma}

  \begin{proof}
We obtain $\theta^*(f)$  simply by replacing the coordinate $x^i$ with $\tilde{x}^i$. Therefore, 
\begin{align*}
w.\theta^*(f)(\tilde{x}^1,\tilde{x}^2,\ldots,\tilde{x}^n)&=\theta^*(f)(\rho_{new}(w)\tilde{x^1},\rho_{new}(w)\tilde{x^2},\ldots,\rho_{new}(w)\tilde{x^n})\\
&=\theta^*(f)\left(\rho_{sign}(w)\rho_{ref}(w)\tilde{x^1},\rho_{sign}(w)\rho_{ref}(w)\tilde{x^2},\ldots,\right.\\
&\left. \ \ \ \ \rho_{sign}(w)\rho_{ref}(w)\tilde{x^n}\right)\\
&=(\rho_{sign}(w))^{deg(f)} \theta^*(f)(\tilde{x}^1,\tilde{x}^2,\ldots,\tilde{x}^n).
\end{align*}
 \end{proof}

Let $z^1,...,z^r$ be algebraically independent  invariant polynomials of $\rho_{new}$ and $u^1,...,u^r$ be the generators of ${\mathbb C}[\rho_{ref}]$ (in the notation of section \ref{first section}). We assume $z^1=\theta^*(u^1)$. Hence, the Hessian of $z^1$ defines a contravariant flat metric $\Omega_2$ on ${\mathcal O}(\rho_{new})$. Examples show that the entries of $\Omega_2(z)$ are  rational  in general and it is hard to construct flat pencil of metrics. We overcome this problem by defining   certain invariants for $\rho_{new}$ which also leads to the construction of Frobenius manifold structures.

\begin{proposition}\label{pb1}
There exist $r$ algebraically independent invariant polynomials $z^1,z^2,\ldots,z^r$ of $\rho_{new}$ with the degrees \be\label{degree new} 2,4,6,\ldots ,2\lfloor \frac{r+1}{2}\rfloor;\  6,8,\ldots,2\lceil \frac{r+3}{2}\rceil.\ee 
\end{proposition}

\begin{proof}
We will use the invariants $u^1,...,u^r$ of $\rho_{ref}$ to construct invariants of $\rho_{new}$. We set $I=\{i:\eta_i \mathrm{~is~ even}\}$ and $J=\{j:\eta_j  \mathrm{~ is ~ odd}\}$. Using Lemma \ref{lem1}, $\theta^*(u^i)$ is an invariant of $\rho_{new}$ for any  $i\in I$. Let  $\kappa$ be the minimal index in $J$. Then $\theta^*(u^\kappa u^j)$ is an invariant of $\rho_{new}$ for any $j\in J$. By this way, we construct $r$ invariants polynomial, $z^1,\ldots,z^r$  for $\rho_{new}$ with the degrees given in  \eqref{degree new}. Note that  any polynomial in $z^1,\ldots,z^r$ can be written as  a polynomial in $u^1,\ldots ,u^r$. Hence,  $z^1,\ldots,z^r$ are algebraically independent.  
\end{proof}

\begin{remark} 
We observe that the invariant polynomials  constructed by Proposition  \ref{pb1} do  not necessarily form a set of primary invariant polynomials of $\rho_{new}$. According to the invariant theory \cite{dresken}, the product of the degrees of  primary invariants is divisible by the order of the group. For example, when $\mathcal W$ is type $A_4$, the degrees of $z^i$ are   $2,4,6,8$. The product of these degrees  is not divisible by the order 120 of the group.
\end{remark} 

We  keep the notations $z^1,...,z^r$ for the invariant polynomials of $\rho_{new}$ constructed in  Proposition  \ref{pb1}.

\begin{theorem}
 The orbits space ${\mathcal O}(\rho_{new})$  has natural  Frobenius manifold structures   isomorphic to the natural Frobenius manifolds structures defined on $ {\mathcal O}(\rho_{ref})$ by Theorem \ref{polyFrob} and Theorem \ref{polyFrob2}.
\end{theorem}
\begin{proof}
We consider the map $(u^1,...,u^r)\to (z^1,z^2,\ldots,z^r)$ given in  Proposition \ref{pb1} as diffeomorphism on some open subset of $u^\kappa\neq 0$ where $\kappa$ is defined in the proof of Proposition \ref{pb1}. Note that, under this diffeomorphism, the metric defined by the Hessian of $u^1$ is identified with the metric defined by the Hessian of $z^1$. Thus, we can transfer  to   ${\mathcal O}(\rho_{new})$, any  regular QFPM  given by the Theorems \ref{polyFrob} and \ref{polyFrob2}. In this way, we obtain natural Frobenius manifold  structures on   ${\mathcal O}(\rho_{new})$.
\end{proof}

\begin{example} \label{A4}
The irreducible reflection representation $\rf$ of Coxeter group of type $A_4$ is generated by the matrices
\be \sigma=\begin{pmatrix}
1&0&0&-1\\
0&1&0&-1\\
0&0&1&-1\\
0&0&0&-1\end{pmatrix}\ \ and \ \  \tau=\begin{pmatrix}
  -1&1&0&0\\
  -1&0&1&0\\
  -1&0&0&1\\
  -1&0&0&0
\end{pmatrix}\ee 
The polynomial ring $\C[\rf]=\C[u_1,u_2,u_3,u_4]$ where  
\begin{eqnarray*}
 u_1 &=&  x_1^2 - \frac{1}{2}x_1x_2 - \frac{1}{2}x_1x_3 -\frac{1}{2}x_1x_4 + x_2^2 -\frac{1}{2}x_2x_3 - \frac{1}{2}x_2x_4 +
        x_3^2 - \frac{1}{2}x_3x_4 + x_4^2,\\ \nonumber
          u_2 &=& x_1^3 - \frac{3}{4}x_1^2x_2 - \frac{3}{4}x_1^2x_3 - \frac{3}{4}x_1^2x_4 - \frac{3}{4}x_1x_2^2 + x_1x_2x_3 +
        x_1x_2x_4 - \frac{3}{4}x_1x_3^2 + x_1x_3x_4 - \frac{3}{4}x_1x_4^2 + x_2^3 - \frac{3}{4}x_2^2x_3 \\ \nonumber
      & & -\frac{3}{4}x_2^2x_4 - \frac{3}{4}x_2x_3^2 + x_2x_3x_4 - \frac{3}{4}x_2x_4^2 + x_3^3 - \frac{3}{4}x_3^2x_4
        - \frac{3}{4}x_3x_4^2 + x_4^3   ,\\ \nonumber
   u_3 &=& x_1^4 - x_1^3x_2 - x_1^3x_3 - x_1^3x_4 + x_1^2x_2x_3 + x_1^2x_2x_4 + x_1^2x_3x_4 -         x_1x_2^3 + x_1x_2^2x_3 + x_1x_2^2x_4 + x_1x_2x_3^2  -3x_1x_2x_3x_4 \\ \nonumber
        & & +  x_1x_2x_4^2 - x_1x_3^3 + x_1x_3^2x_4 + x_1x_3x_4^2 - x_1x_4^3 + x_2^4 -
        x_2^3x_3 - x_2^3x_4 + x_2^2x_3x_4 - x_2x_3^3 + x_2x_3^2x_4+ x_2x_3x_4^2    \\ \nonumber
        & &  -   x_2x_4^3 + x_3^4 - x_3^3x_4 - x_3x_4^3 + x_4^4, \\ \nonumber
   u_4 & = &  x_1^5 - \frac{5}{4}x_1^4x_2 - \frac{5}{4}x_1^4x_3 - \frac{5}{4}x_1^4x_4 + \frac{5}{3}x_1^3x_2x_3 +
        \frac{5}{3}x_1^3x_2x_4 + \frac{5}{3}x_1^3x_3x_4 - \frac{5}{2}x_1^2x_2x_3x_4 - \frac{5}{4}x_1x_2^4      + \frac{5}{3}x_1x_2^3x_3 \\ \nonumber
        & &  + \frac{5}{3}x_1x_2^3x_4 - \frac{5}{2}x_1x_2^2x_3x_4 + \frac{5}{3}x_1x_2x_3^3 -
        \frac{5}{2}x_1x_2x_3^2x_4 - \frac{5}{2}x_1x_2x_3x_4^2 + \frac{5}{3}x_1x_2x_4^3 - \frac{5}{4}x_1x_3^4 +
        \frac{5}{3}x_1x_3^3x_4  \\ \nonumber
        & &  + \frac{5}{3}x_1x_3x_4^3 - \frac{5}{4}x_1x_4^4 + x_2^5   - \frac{5}{4}x_2^4x_3     -  \frac{5}{4}x_2^4x_4 + \frac{5}{3}x_2^3x_3x_4 - \frac{5}{4}x_2x_3^4 + \frac{5}{3}x_2x_3^3x_4 +
       \frac{5}{3}x_2x_3x_4^3 -\frac{5}{4}x_2x_4^4 + x_3^5 \\ \nonumber
        & & - \frac{5}{4}x_3^4x_4  - \frac{5}{4}x_3x_4^4 + x_4^5.
\end{eqnarray*}

The Frobenius manifold of type $A_4$ is a result of the regular QFPM consists of $\Omega_2(u)$  and  $\Omega_1={\partial_{u_4}} \Omega_2(u)$ where $\Omega_2(u)$ is defined by the Hessian of $u_1$.  The representation   $\rho_{new}$ is generated by $\tau$ and $-\sigma$. Then the primary invariants of  $\rho_{new}$ have degrees $2,4,6,10$ while the secondary invariants have  degrees $8,13,15$.  The Hessian of the degree 2 invariant  $z_1$ leads to the flat contravariant  metric $ \Omega_2(z)$ but it is hard to find a FPM. We fix the following $4$ invariants polynomials for  ${\mathcal O}(\rho_{new})$ of degrees $2,4,6$ and $8$: \[{z_1}=u_1,\ \ {z_2}=u_3,\ \ {z_3}=u_2^2, \ \ z_4= u_2 u_4.\]
Then the matrix of $\Omega_2(z)$ consists of the columns  
$${\Omega_2^{i1}(z)}=\begin{pmatrix}
z_1		\\
2 z_2\\
3 z_3\\
4 z_4\end{pmatrix}, \ \ {\Omega_2^{i2}(z)}=\begin{pmatrix}
 2 z_2	\\
	-\frac{64}{625} z_1^3+\frac{68}{25} z_1 z_2+\frac{864}{625}z_3\\
	\frac{12}{5} z_1 z_3+\frac{18 }{5}z_4\\
	\frac{64}{75} z_1^2 z_3+\frac{43}{15} z_2 z_3+\frac{62}{25} z_1 z_4+\frac{9 }{5 }\frac{ z_4^2}{ z_3}
\end{pmatrix}, \ \ $$ $$ {\Omega_2^{i3}(z)}=\begin{pmatrix}
 3 z_3	\\
 	\frac{12}{5} z_1 z_3+\frac{18 }{5}z_4\\
\frac{2}{3} z_1^2 z_3+\frac{25}{3} z_2 z_3\\
	-\frac{26}{45} z_1^3 z_3+\frac{95}{18} z_1 z_2 z_3+\frac{14 }{5}z_3^2+\frac{1}{3} z_1^2 z_4+\frac{25}{6} z_2 z_4
\end{pmatrix}$$
and
$${\Omega_2^{i4}(z)}=\begin{pmatrix}
 4 z_4\\
	\frac{64}{75} z_1^2 z_3+\frac{43}{15} z_2 z_3+\frac{62}{25} z_1 z_4+\frac{9 }{5 }\frac{ z_4^2}{ z_3}\\
-\frac{26}{45} z_1^3 z_3+\frac{95}{18} z_1 z_2 z_3+\frac{14 }{5}z_3^2+\frac{1}{3} z_1^2 z_4+\frac{25}{6} z_2 z_4\\
	\frac{214 }{2025}z_1^4 z_3+\frac{52}{81} z_1^2 z_2 z_3+\frac{625}{324} z_2^2 z_3+\frac{56}{25} z_1 z_3^2-\frac{26}{45} z_1^3 z_4+\frac{95}{18} z_1 z_2 z_4+\frac{62}{15} z_3 z_4+\frac{1}{ z_3}\frac{z_1^2 z_4^2}{ z_3}+\frac{25 }{12 }\frac{ z_2 z_4^2}{ z_3}
\end{pmatrix}$$
Therefore, on ${\mathcal O}(\rho_{new})$, we get the regular QFPM formed by $\Omega_2(z)$ and $ {\Omega}_1(z)=\L_{e} {\Omega_2}$ where $e=\sqrt{z_3}\partial_{z_4}$.  Of course, the resulted Frobenius manifold is  of type $A_4$.
\end{example}

\begin{remark} 
It is straightforward to  generalized the results of this section to other types of Coxeter groups and we obtain natural Frobenius manifolds on  ${\mathcal O}(\rho_{new})$. But we  lack sorting out when $\rho_{new}$ is not a reflections group (i.e. see Proposition \ref{new not ref}). Robert Howett informed us that when $\mathcal W$ is of type $E_8$, the representation $\rho_{new}$ is generated by reflections.
\end{remark}

\section{Dihedral and dicyclic  groups}\label{dicyclic and dihedral} 
In this section, we give results of applying Dubrovin's method to irreducible representations of the dihedral groups  (Coxeter groups of type) $I_2(m)$, $m>2$ and Dicyclic groups $Dic_m$.
 We mention that Dubrovin computed by an ad-hoc procedure  all possible potentials of 2-dimensional Frobenius manifolds  \cite{DuRev}. Here we find some of them are related to invariant theory of finite groups.

\subsection{Dihedral groups}\label{dihedral}

Irreducible representations of $I_2(m)$ are of  rank 1 or 2. Let   $\xi_m$ be a primitive $m$-th root of unity. The rank 2 representations are $\rho_k$, $k=1,2,\ldots,\frac{m-2}{2}$ for even $m$, and $k=1,2,\ldots,\frac{m-1}{2}$ for odd $m$. Here,  $\rho_k$ is generated by the matrices 
\begin{equation}\
 \begin{pmatrix} \xi_m^k &0\\
0&\xi_m^{-k}\end{pmatrix}
,~ \begin{pmatrix} 0&1\\
1& 0\end{pmatrix}.
\end{equation}
When $k=1$, we get the standard reflection representation of $I_2(m)$.
 We observe that    $\mathbb{C}[\rho_k]$ can be interpreted as the invariant ring of the standard reflection representation of $I_2(h)$  where $h=\frac{m}{\gcd(m,k)}$, i.e. it is generated by
 \[t_1=\frac{1}{h} x_1 x_2,\ t_2=x_1^h+x_2^h.\]
 
 Hence, applying Dubrovin's method, we  get the polynomial Frobenius manifold of type $I_2(h)$  and its conjugate $\widetilde{I}_2(h)$ obtained by Theorem \ref{main thm}.

\subsection{Dicyclic groups}\label{dicyclic}

We fix  a natural number $m$. The dicyclic group $\mathrm{Dic}_m$ is a group of order $4m$ defined by
 \begin{equation}
 \mathrm{Dic}_m=\langle\sigma,\alpha|\sigma^{2m}=1, \alpha^2=\sigma^m, \alpha^{-1} \sigma \alpha=\sigma^{-1}\rangle. 
 \end{equation}

 The irreducible representation of $\mathrm{Dic}_m$ are of rank 1 or 2.  The 2-dimensional irreducible representations  $\psi_k$ and $\varrho_{l}$ are defined by setting \begin{equation}\label{dic1}
\psi_k(\sigma)=\left(\begin{array}{cc}
     \xi_{2m}^k& 0  \\
0     &\xi_{2m}^{-k}
\end{array}\right),~\psi_k(\alpha)=\left(\begin{array}{cc}
    0 & 1 \\
    -1 & 0
\end{array}\right),
\end{equation}
and 
\begin{equation}\label{dih1}
\varrho_{l}(\sigma)=\left(\begin{array}{cc}
     \xi_m^k& 0  \\
0     & \xi_m^{-k}
\end{array}\right),~\varrho_{l}(\alpha)=\left(\begin{array}{cc}
    0 & 1 \\
    1 & 0
\end{array}\right).
\end{equation}
Here $1\leq k\leq \frac{m-2}{2}$ and $1\leq l\leq m-1$ when $m$ is even while $1\leq k\leq \frac{m-1}{2}$ and $1\leq l\leq m-2$ when $m$ is odd.
Note that $\psi_1$ is the standard representation of $\mathrm{Dic}_m$ in the litreture. We observe that the invariant ring ${\mathbb C}[\varrho_l]$ can be interpreted as the invariant ring $\C[\rho_l]$ where $\rho_l$ is the representation of $I_2(m)$ given in  section \ref{dihedral}. Thus, the result of applying Dubrovin's method to  $\varrho_l$ is given in that section.  We consider here  the representations $\psi_k$. Let us fix the integer $k$ and  set $h={m\over \gcd(m,k)}$. We define 
\begin{equation} \label{dicin} u_1=x_1^2x_2^2,~ u_2= x_1^{2h}+x_2^{2h},~u_3=x_1x_2(x_1^{2h}-x_2^{2h}).
\end{equation} 
It is straightforward to verify  that $u_1$, $u_2$ and $u_3$ are invariants under the action of $\psi_k$.
\begin{proposition}\label{dihinv}
The invariant ring ${\mathbb C}[\psi_k]$ is generated by $u_1,u_2$ and $u_3$.
\end{proposition}
\begin{proof}
 A general homogeneous polynomial of degree $q$ has the form  \begin{equation} f(x_1,x_2)=a_q x_1^q + a_{q-1} x_1^{q-1} x_2 + \cdots + a_1 x_1 x_2^{q-1} + a_0 x_2^q\end{equation}  where  $a_0 , \cdots , a_q \in {\mathbb C} $. Being  invariant under $\psi_k(\alpha)$, we get 
\begin{align*}
f&=a_q x_1^q + a_{q-1} x_1^{q-1} x_2 + \cdots + a_1 x_1 x_2^{q-1} + a_0 x_2^q\\
=\psi_k(\alpha)f&=(-1)^q a_q x_2^q +(-1)^{q-1} a_{q-1} x_2^{q-1} x_1 + \cdots - a_1 x_2 x_1^{q-1} + a_0 x_1^q
\end{align*}
Thus $a_i = (-1)^{q-i} a_{q-i}$  for  all $i=0 , \cdots ,q$. Similarly, the invariant of  $f$  under $\psi_k(\alpha^2)$ implies   $q$ is even. Hence, $f$ has the form 
\begin{align*}
f&=\sum_{i=0}^{ \frac{q}{2}} a_{q-i} (x_1 x_2)^i [x_1^{q-2i} + (-1)^{q-i} x_2^{q-2i}].
\end{align*}
Moreover,
\begin{align*}
f&=\psi_k(\sigma)f=\sum_{i=0}^{ \frac{q}{2}} a_{q-i}(x_1 x_2)^i [\xi^{-k(q-2i)} x_1^{q-2i} +  (-1)^{q-i}\xi^{k(q-2i)} x_2^{q-2i}]
\end{align*}
implies  $k(q-2i) = 0\ mod(2m)$. Then  $q-2i=2hl$ for some integer $l$.
Therefore, we can write
 \begin{equation}
f=\sum_{q=2hl+2i} a_{q-i}(x_1 x_2)^{i} [ x_1^{2hl} + (-1)^{q-i}  x_2^{2hl}].
\end{equation}
Now we show that  $f\in \mathbb{F}[u_1,u_2,u_3]$. It is sufficient to prove  $\widetilde{f}_l= x_1^{2hl} +   x_2^{2hl}$ and $ \widehat{f}_l= x_1 x_2(x_1^{2hl} -   x_2^{2hl})$ are invariant for every natural number $l$.  When  $l=1$, $\widetilde{f}_l=u_2$ and $\widehat{f}_l=u_3$. For $l+1$, we get
\begin{align}
\widetilde{f}_{l+1}&= x_1^{2h(l+1)} + x_2^{2h(l+1)}= (x_1^{2h} + x_2^{2h})^{l+1}- \displaystyle{\sum_{d=1}^{l}} \left( ^{l+1}_{\ \ d} \right)
x_1^{2hd} x_2^{(l+1-d)2h}\\\nonumber
&= (x_1^{2h} + x_2^{2h})^{l+1}-  \displaystyle{\sum_{d=1}^{\lfloor  \frac{l}{2} \rfloor} } \left( ^{l+1}_{\ \ d} \right)  (x_1 x_2 )^{2hd} ( x_2^{2h(l+1-2d)}+  x_1^{2h(l+1-2d)} ).
\end{align}
Since $d \geq 1$, we have $l+1 -2d \leq l-1 <l $. Therefore, by the induction   $\widetilde{f}_{l+1}\in {\mathbb C}[u_1,u_2,u_3]$. Likewise   $\widehat{f}_{l+1}\in {\mathbb C}[u_1,u_2,u_3]$ since
\begin{align}
\widehat{f}_{l+1}&=x_1x_2(x_1^{2h}-x_2^{2h})(x_1^{2hl}+x_1^{2h(l-1)}x_2^{2h}+x_1^{2h(l-2)}x_2^{4h}+\ldots +x_2^{2hl})\\\nonumber
&=x_1x_2(x_1^{2h}-x_2^{2h})[(x_1^{2hl}+x_2^{2hl})+(x_1x_2)^{2h}(x_1^{2h(l-2)}+x_2^{2h(l-2)})\\ \nonumber
& + (x_1x_2)^{4h}(x_1^{2h(l-4)}+x_2^{2h(l-4)})+\ldots].
\end{align}
This proves the proposition.
\end{proof}

We note that the invariant ring ${\mathbb C}[\psi_k]$ can be interpreted as the invariant ring of  the standard representation of $\mathrm{Dic}_h$. A result of applying Dubrovin's method  is obtained in \cite{dicy}. We summarize the construction here. 

The flat contravariant metric defined by the inverse of the Hessian of $u_1$ is 
\begin{equation}
\Omega_2(u)=\left(
\begin{array}{cc}
 \frac{4 }{3}u_1 & \frac{2h}{3}  u_2 \\
 \frac{2h}{3}  u_2 & -\frac{2 h^2}{3u_1} (  u_2^2-6 u_1^h)
\end{array}
\right).
\end{equation}
Then we considered a  vector field $e$ in the form  $e=f(u_1)\partial_{u_2}$ and imposed the conditions  $\Lie_e\Omega_2$ is flat and $\Lie_e^2\Omega_2=0$. These conditions lead to two independent solutions
\begin{equation}
f_{\pm}=u_1^{\frac{h}{2} \left(1\pm\sqrt{3}\right)}.
\end{equation}
Setting  $e_\pm=f_{\pm}\partial_{u_2}=u_1^{\frac{h}{2} \left(1\pm\sqrt{3}\right)}\partial_{u_2}$, we get  regular quasihomogenous flat pencils of metrics  $(\Omega_2,\Lie_{e_\pm}\Omega_2)$ of degree $
d=\frac{\sqrt{3}h \pm 2}{\sqrt{3}h}$ with $\tau=\mp{\sqrt{3}\over 2h}u_1$. The resulting Frobenius manifold structures are conjugate to each other. The corresponding flat coordinates of reads
\begin{equation} t_1=\mp\frac{\sqrt{3}}{2 h} u_1,~~t_2=u_2 u_1^{\frac{\mp h}{2} \left(\sqrt{3}\pm 1\right)
   }
\end{equation}
Which lead to the potentials  \begin{equation}\label{posPotential}
{\mathbb{F}}=\frac{2^{\mp\sqrt{3} h}
   3^{\frac{1}{2}
   \left(1\pm\sqrt{3} h\right)}
   \left(h
   t_1\right){}^{1\mp\sqrt{3}
   h}}{\mp(3 h^2-1)}+\frac{1}{2}
   t_1 t_2^2\end{equation}
   of the degrees $ \mp\frac{2}{\sqrt{3} n}$ and 1.

\subsection{Finite subgroups of $SL_2({\mathbb C})$}\label{sl2}
In this section we use Dubrovin's method on finite non trivial subgroups of $SL_2({\mathbb C})$. They are classified up to conjugation  and they are called binary polyhedral groups. They consist of the cyclic groups $\mathcal C_m$ and  binary dihedral groups $\mathcal D_m$, binary tetrahedral group $\mathcal T$, binary octahedral group $\mathcal O$ and binary icosahedral group $\mathcal{I}$. We treat them as representations of the corresponding groups. It is known that the invariant rings of these representations are not  polynomial rings and the relations between the generators lead to the classification of simple hypersurface  singularities. We use the  sets of generators of the invariant rings listed in \cite{Leuschke}. Applying Dubrovin's method, we obtain natural polynomial Frobenius manifold structure and their conjugations (as given in section \ref{first section}).  We write  below  only  the flat coordinates and the type of the resulting polynomial Frobenius manifold structures.  Note that the findings are not apparent from examining the invariant rings.
\begin{enumerate}
\item Cyclic groups $\mathcal C_m$: Here $m\geq 2$ and the invarinat ring is generated by $xy,~x^m,~y^m$. We fix the following invariant polynomials  
 \[t_1=\frac{1}{m} x  y,\ t_2=x^m+y^m\]
Then the ring generated by $t_1$ and $t_2$ is isomorphic to the invariant ring of the standard representation of the dihedral group $I_2(m)$. Thus, using Dubrovin's method and $(t_1,t_2)$ as coordinates on the orbits space, we get Frobenius manifold of type  $I_2(m)$. 

In case we set  $t_1=\frac{1}{m} x  y$ and $t_2=x^m$, we get the WDVV solution  $\frac{1}{2} t_1 t_2^2$. It corresponds to a trivial Frobenius manifold structure but here the natural charge is  $\frac{m-2}{m}$ while the degrees are $\frac{1}{m}$  and 1. 

\item {The binary dihedral group  $\mathcal D_m$:} This is the standard representation of the dicyclic group  $\mathrm{Dic}_m$. A result of applying Dubrovin's method is given in section \ref{dicyclic}.
\item {The binary tetrahedral $\mathcal T$:} 
We fix the following set of generators of the invariant ring
 \begin{equation} t_1={5\over 12}x y \left(x^4-y^4\right), \ \ t_2=\left(x^4+y^4\right)^3-36 x^4 y^4 \left(x^4+y^4\right).\end{equation}
 \begin{equation*}  t_3=16 x^4 y^4+2 \left(x^4-y^4\right).\end{equation*}
We choose $(t_1,t_2)$ as coordinates on the orbits space. Then  the Hessian of ${12\over 5}t_1$ defines a flat metric $\Omega_2(t)$ linear in $t_2$. Here, we apply Lemma \ref{almost linear} and we get  a regular QFPM of degree $d=\frac{1}{2}$ with $\tau= t_1$ consists of 
\begin{equation}   \Omega_2^{ij}(t)=\left(
\begin{array}{cc}
 \frac{1}{2}t_1 & t_2 \\
 t_2 & \frac{-4478976}{625}  t_1^3 \\
\end{array}
\right),\Omega_1^{ij} = \Lie_{\partial_{t_2}} \Omega_2^{ij}(t)=\left(
\begin{array}{cc}
 0 & 1 \\
 1 & 0 \\
\end{array}
\right) \end{equation} 
The resulting Frobenius manifold is  of type $I_2(4)$. 
\item {The binary octahedral $\mathcal O$:}  Let us fix the generators  of the invariant ring to be \begin{equation} t_1={7\over 12}(16 x^4 y^4+\left(x^4-y^4\right)^2),\ \  t_2=\left(x y \left(x^4-y^4\right)\right)^2, \end{equation}
\begin{equation*}t_3=y x^{17}-34 y^5 x^{13}+34 y^{13} x^5-y^{17} x\end{equation*}
In the coordinates $(t_1,t_2)$,  the metric $\Omega_2(t)$ defined by the Hessian of ${12\over 7} t_1$ is linear in $t_2$ and leads to a regular QFPM  with charge $\frac{1}{3}$.   The resulting Frobenius manifold  is  of type $I_2(3)$. 

\item {The binary icosahedral $\mathcal I$:} We fix the generators of the  invariant ring  \begin{align}
t_1&={11\over 30}(x^{11} y+11 x^6 y^6-x y^{11}),\\\nonumber
t_2&=x^{30}+522 x^{25} y^5-10005 x^{20} y^{10}-10005 x^{10} y^{20}-522 x^5 y^{25}+y^{30},\\\nonumber
t_3 &=x^{20}-228 x^{15} y^5+494 x^{10} y^{10}+228 x^5 y^{15}+y^{20} 
\end{align}
Fix $(t_1,t_2)$ as coordinates, the Hessian of ${30\over 11}t_1$ leads to a metric $\Omega_2(t)$ linear in $t_2$. The regular QFPM formed by  $\Omega_2$ and $\Omega_1 = \partial_{t_2}\Omega_2(t)$ leads to  a  Frobenius manifold  of type $I_2(5)$. 

\end{enumerate}

  \section{Finite subgroups of $SL_3({\mathbb C})$} \label{finite rank 3}
 Finite subgroups  of $SL_3({\mathbb C})$ are classified into the families  $(\mathcal A)$, $(\mathcal B)$, \ldots, $(\mathcal L)$ \cite{yau}. We treat them as representations of the corresponding groups and they are not reflection representations.  Watanabe and  Rotillon listed in \cite{watan}  those subgroups   where the invariant rings  are complete intersections missing  type $(\mathcal J)$ and $(\mathcal K)$. These missing groups were recognized by  Yau and Yu  \cite{yau}. In the end, there is a total of 29 types of  finite subgroups of $SL_3({\mathbb C})$ whose invariant rings are complete intersection and  their sets of generators are known explicitly. We treat them as linear representations of fnite groups and we apply Dubrovin's method. The set of generators is taken from \cite{watan} and we use the same numbering  $(1)$, $(2)$, \ldots, $(27)$  of the 27 families of subgroups listed there.
 
 Recall that to apply Dubrovin's method, we must find 
 \begin{equation} \begin{array}{c}\label{cond}
\text{ a minimal degree invariant polynomial where the Hessian defines}\\
 \text{a flat contravariant metric.}
 \end{array}\end{equation} 
 This condition excluded the following subgroups
 \begin{enumerate}
 \item $(17)$ which are of type $(\mathcal A)$.
 \item $(3)-(8)$, $(19)-(23)$ which are of type $(\mathcal B)$.
 \item $(10)$, $(24)$ and $(25)$ which are of type $(\mathcal C)$.
 \item $(13)$, $(15)$, $(16)$ and $(27)$ which are of types $(\mathcal G)$, $(\mathcal{L})$, $(\textbf{I})$ and $( \mathcal{ E})$, respectively. 
 \item The groups $(\mathcal J)$ and $(\mathcal K)$ which are not considered in  \cite{watan}.
 \end{enumerate}
 
 For the remaining family of subgroups, when condition \eqref{cond} is satisfied,  we use  Lemma \ref{almost linear} to construct flat pencil of metric under   appropriate choice  of a set of invariant polynomials. We  find  natural Frobenius manifold structures of types $A_3$, $B_3$, $H_3$, $B_3^1$, or the trivial $T_3$. In each case,  we will mention  the type of the resulting Frobenius structure and the corresponding flat coordinates.  From Theorem  \ref{conj on orbits}, we know that ones the orbits space acquire one of these structures then it also possess the conjugate structure (see Example \ref{Dualtiy poly} and  Example \ref{Triv1}).   Thus we will not mention explicitly the appearance of the natural conjugate Frobenius manifold  structures.  
 
 \begin{enumerate}
     \item[$(1)$] This is a family of  groups  of type  $(\mathcal A)$ depending  on an integer $m>1$. Complete set of generators of the invariant rings  consists of  $x^m,y^m,z^m,xyz$. The Hessian of $xyz$ does not define a flat metric. Hence, condition \eqref{cond} exclude the case $m>3$.
     
 For $m=2$, we fix the invariant polynomials 
\begin{equation} u_1=x^2+y^2+z^2,\ u_2=x^2 y^2+z^2 y^2+x^2 z^2,\ u_3=(x y z)^2.
\end{equation}
Then $\{u_1,u_2,u_3\}$ can be identified with the set of generators of the invariant  ring of the standard reflection representation of Coxeter groups of type $B_3$. Thus, applying Dubrovin's method, we get  natural Frobenius  manifold structures of  types $A_3$, $B_3$ and $B_3^1$. We also get the natural  trivial Frobenius manifold structure of type $T_3$ using the setting of the family  (2) given below. Thus, considering the conjugate structures and Frobenius manifolds obtained in Example \ref{Triv1}, we proved that the orbits space has 8 different natural Frobenius manifold  structures. 

For $m=3$, we fix the invariant polynomials 
\begin{equation} u_1=x^3 + y^3 + z^3,\ u_2 = x^3 y^3 + y^3 z^3 +z^3 x^3,\ u_3 = (x y z)^3.
\end{equation} 
The Hessian of $u_1$ defines a contravariant flat metric $\Omega_2$. This metric and its Christoffel symbols are almost linear in each variable $u_i$. We can and will apply Lemma \ref{almost linear} and we get three regular QFPM. 
From QFPM  $(\Omega_2,\L_{\partial_{u_3}}\Omega_2)$, we get Frobenius manifold structure of type $B_3$. It has  flat coordinates  \begin{equation} t_1=\frac{2 }{9}u_1,\ t_2=-\frac{u_1^2-4 u_2}{6 \sqrt{2}},\ t_3=\frac{7 u_1^3}{216}-\frac{1}{6} u_2 u_1+u_3.\end{equation}
The QFPM $(\Omega_2,\L_{\partial_{u_2}}\Omega_2)$ leads to type $A_3$. It has  flat coordinates 
    \begin{equation} 
    t_1={1\over 3} u_1,\ t_2=u_2-{1\over 8} u_1^2,\ t_3=\sqrt{u_3}.
    \end{equation}   
Finally we get Frobenius manifold structure of type $B_3^1$ from the QFPM  $(\Omega_2,\L_{\partial_{u_1}}\Omega_2)$. Here the flat coordinates are  \begin{equation}\label{rational 1 m=3 }  t_1=u_1,\ t_2=u_2 u_3^{-{1\over 4}}, t_3={4\over 3}u_3^{1\over 4}.\end{equation}

\item[$(2)$] This is a family of groups of type $(\mathcal B)$ depending on an integer  $m \geq 1$. The polynomials $x^{2 m}+y^{2 m},\  (x y)^2$, $x y z \left(x^{2 m}-y^{2 m}\right)$ and $z^2$ form  complete  sets of generators for the invariant rings. Because of condition \eqref{cond}, we need only to consider  $m =1$. In this case,    we fix the invariant polynomials
\[u_1= x^2+y^2+z^2,
u_2=z^2,
u_3 =x^2 y^2
\]
The metric $\Omega_2(u)$ defined by the Hessian of $u_1$ and its Christoffel symbols are linear in each variable $u_i$. However, Lemma \ref{almost linear} is applicable only for $u_2$. The QFPM $(\Omega_2,\L_{\partial_{u_2}} \Omega_2)$ has degree 0 with $\tau=u_1$. It leads to a natural trivial Frobenius manifold structure of type $T_3$. Here the 
 flat coordinates are 
\begin{equation} 
t_1={1\over 2} u_1,\ t_2=u_2-{1\over 2} u_1,\ t_3= (-2u_3)^{1\over 2}.
\end{equation} 
  \item[$(9)$] This is a family of groups  of type $(\mathcal C)$ depending on an integer  $m> 1$. Complete sets of generators of the invariant rings consist
of  $ x y z,x^m+y^m+z^m,x^m y^m+x^m z^m+y^m z^m$, and $(x^m-y^m)(z^m-x^m) (y^m-z^m)$. Here we get the same natural Frobenius manifold structure obtained for the family (1).
 
\item[$(11)$] This is family of  groups  of type $(\mathcal C)$ depending on an integer $m>1$.  Complete  sets of generators of the  invariant rings consists of  
\begin{equation} 
u_1= x^m+y^m+z^m,\ u_2=x^2 y^2 z^2, \ u_3= x^m y^m +y^mz^m+z^m x^m  
 \end{equation} 
 and 
 \begin{equation} 
 u_4=x y z \left(x^m-y^m\right) \left(z^m-x^m\right)
   \left(y^m-z^m\right).
 \end{equation} 
Since the Hessian $u_2$ does not define a flat metric, we consider only $2\leq m\leq 6$. For $m=2$ we can use the same argument given for the family (1).

For $3\leq m\leq 6$, the contravariant metric $\Omega^{ij}_2$ defined by the Hessian of $u_1$  and its Christofel symbols are almost linear in $u_1$ and $u_3$ and we can apply Lemma \ref{almost linear} to both variables.

The FPM  $(\Omega_2,\L_{\partial_{u_3}}\Omega_2)$ is regular  quasihomogeneous of degree $1\over 2$ with $\tau=u_1$. We can fix the flat coordinates 
     \begin{equation} 
     t_1={m-1\over 2m} u_1,\ t_2=u_2^{m\over 4},\ t_3= u_3-{1\over 8} u_1^2.
     \end{equation} 
The resulting natural Frobenius manifold structure is a  polynomial  of type $A_3$.  

Similarly, the FPM  $(\Omega_2,\L_{\partial_{u_1}}\Omega_2)$ is regular  quasihomogeneous of degree $0$ with $\tau={m-1\over m} u_1$. We can fix the flat coordinates  
\begin{equation}\label{rational 11}
t_1=u_1,\ t_2=u_2^{m\over 8},\ t_3=u_3 u_2^{-{m\over 8}}.
\end{equation} 
We get a natural Frobenius manifold structure  of type $B_3^1$. 

\item[$(12)$] This is a group of type  $(\mathcal F)$.  A  complete  set of generators of the invariant ring consists of  \begin{equation} u_1=(x^3 + y^3 + z^3)^2-12(x^3 y^3 + y^3 z^3 + z^3 x^3),u_2= (x^3 - y^3) (y^3 - z^3) (z^3 - x^3),\end{equation}
 \begin{equation*}
 u_3=( x y z)^4+216 (x y z)^3 (x^3 + y^3 + z^3), u_4= ((x^3 + y^3 + z^3)^2 - 18 (x y z)^2)^2. \end{equation*} 
The FPM  $(\Omega^{ij}_2,\L_{\partial_{u_3}} \Omega^{ij}_2)$ is  regular quasihomogeneous of degree $d={1\over 2}$ with $\tau= {5\over 12} u_1$. Flat coordinates are given by \begin{equation} t_1=\frac{5 }{12}u_1,\ t_2=10  \sqrt{\frac{62}{41}} u_2,\ t_3=u_3-\frac{847 }{1312}u_1^2.\end{equation} 
 The resulting natural Frobenius manifold structure is of type $A_3$.

\item[$(14)$] This is a group of type  $(\mathcal H)$ and  a minimal set  of generators of the invariant ring consists of $$u_1=x^2+y z,\ \ u_2=8 y z x^4-2 y^2 z^2 x^2-\left(y^5+z^5\right) x+y^3 z^3,$$
and 
\begin{align*} u_3&=y^{10}+6 z^5 y^5+20 x^2 z^4 y^4-160 x^4 z^3 y^3+320 x^6 z^2 y^2+z^{10}\\&-4 x \left(y^5+z^5\right) \left(32 x^4-20 y z x^2+5 y^2 z^2\right).\end{align*}

The Hessian of $10 u_1$ leads to a regular  QFPM  $(\Omega_2,\L_{ \partial_{u_3}}\Omega_2)$ of degree $d={4\over5}$ with $\tau= {1\over 10} u_1$. By fixing the flat coordinates $$t_1=\frac{1}{10}u_1,\ t_2=\sqrt{2} u_2-\sqrt{2} u_1^3,\ t_3=14 u_1^5-20 u_2 u_1^2+u_3,$$ 
we arrive to a natural polynomial Frobenius structure of type $H_3$.

\item[$(18)$] This is a family of groups of type $(\mathcal B)$ depending on integers $ p\geq 1$ and $ q\geq 2$. A complete  sets of  generators of the invariant rings consists of $(x^{2 p q} + 
  y^{2 p q}), (x y  )^{2 q}, (x y z)^2, (x^{2 p q} -    y^{2 p q}) x y z , z^{2 q}$. The Hessian of $(xyz)^2$ does not define a flat metric.  From  condition \eqref{cond}, we consider only the two cases: $p=1$ but  $q=2$ or  $q=3$. In both cases we get 3 types of natural Frobenius manifold  structures.
  The   first natural Frobenius structure is of type $T_3$. It has the flat coordinates 
  \begin{equation} t_1=\frac{ 2 q-1}{2 q}  (x^{2q}+y^{2q}+z^{2q}),\ t_2=\frac{-1}{2}  (x^{2q}+y^{2q}),\ t_3= (\frac{2-4q}{q}x y)^{\frac{q}{2}}.\end{equation}
 The corresponding  
  regular QFPM is $(\Omega_2, \L_{\partial_{t_2}}\Omega_2)$ with $\tau=  t_1$ where $\O_2$ defined by the Hessian of $t_1$. Let us fix 
  \be u_1=x^{2q}+y^{2q}+z^{2q},\ u_2=  (x y z)^2,\ u_3= -2x^{2q}  y^{2q} -2x^{2q} z^{2q} -2y^{2q} z^{2q}.\ee  Then the second  natural Frobenius manifold structure is of type  $B_3^1$. It has the flat coordinates
  \be t_1=\frac{2q-1}{2q} u_1,\ t_2=u_2^{\frac{q}{4}},\ t_3= u_3u_2^{\frac{-q}{4}}
   \ee   
  The  corresponding  regular QFPM is $(\Omega_2^1,\L_{ \partial_{t_1}}\Omega_2)$ has degree 0  with $\tau=t_1$.
   Finally, we get natural Frobenius manifold structure of type $A_3$ having the flat coordinates 
    \begin{equation} t_1=\frac{2q-1}{4q} u_1,\ t_2=\frac{2 \sqrt{2q-1}}{\sqrt{q}} u_2^{\frac{q}{2}},\ t_3= u_3 +\frac{1}{4} u_1^2.\end{equation} 
   The corresponding  regular  QFPM  $(\Omega_2,\L_{ \partial_{t_3}}\Omega_2)$ is of degree $\frac{1}{2}$ with $\tau=t_1$.

\item[$(26)$] This is a family of  groups of type $(\mathcal C)$ depending on  even integer $m\geq 2$. The set of  minimal generators of invariant ring has  \[ x^{3 m} + y^{3 m} + z^{3 m}, (x y z)^2, x^{2 m} y^m+ x^m y^{2 m} + y^{2 m} z^m + y^m z^{2 m} + z^{2 m} x^m + 
 z^m x^{2 m},\] \[x y z (x^m - y^m) (y^m - z^m) (z^m - x^m),(x^m - y^m)^2 (y^m - z^m)^2 (z^m - x^m)^2.\]
 The only possible case under  condition \eqref{cond} is when $m=2$. In this case we get a natural Frobenius manifold of type $A_3$. It has the flat coordinates  $$t_1=\frac{5 }{12}(x^{6}+y^{6}+z^{6}),\ t_2= \sqrt{\frac{20}{3}} (x y z)^3,\ t_3=x^{12}+y^{12}+z^{12}-\frac{3 }{4}(x^{6}+y^{6}+z^{6})^2.$$
Here, $\O_2$ is defined by the Hessian of ${12\over 5}t_1$ and the corresponding regular QFPM $(\Omega_2, \L_{\partial_{t_3}} \Omega_2)$ is of degree $d={1\over2}$ with $\tau= t_1$.

On the other hand, the Hessian of ${5\over 6}t_1$ leads to a regular QFPM $(\Omega_2, \L_{\partial_{t_1}} \Omega_2)$ of degree 0 with $\tau= t_1$. In this case  the flat coordinates are 
 \begin{equation} t_1=\frac{5}{6} (x^{6}+y^{6}+z^{6}),\ t_2=  (x y z)^{\frac{3}{2}},\end{equation} \begin{equation*} \ t_3= \frac{-5}{6}(x y z)^{\frac{-1}{2}} \left(x^{12}+y^{12}+z^{12}-\frac{3 }{4}(x^{6}+y^{6}+z^{6})^2\right).\end{equation*} 
 The resulting natural Frobenius manifold structure is of type $B_3^1$.

 \end{enumerate}
 

\section*{Acknowledgements}
The  authors thank Robert Howett, Hans-Christian Herbig and Christopher Seaton for useful  discussions.  They very much appreciate the Magma program's support team for their helpful cooperation.

\section*{Funding} This work was partially funded by the internal grant of Sultan Qaboos University (IG/SCI/DOMS/19/08).

\section*{Data Availability} Non applicable. 


\noindent Yassir Dinar \\
\noindent dinar@squ.edu.om \\

\noindent Zainab Al-Maamari\\
\noindent s100108@student.squ.edu.om\\

\noindent Depatment of Mathematics\\
\noindent College of Science\\
\noindent Sultan Qaboos University\\
\noindent Muscat, Oman

\end{document}